\documentclass[preprint,12pt]{elsarticle}

\usepackage{amssymb,xspace}
\usepackage{amsfonts,amsmath,latexsym}
\ifx\pdftexversion\undefined
  \usepackage[a4paper,colorlinks,linkcolor=black,filecolor=black,citecolor=black,urlcolor=black,pdfstartview=FitH]{hyperref}
\else
  \usepackage[a4paper,colorlinks,linkcolor=blue,filecolor=blue,citecolor=blue,urlcolor=blue,pdfstartview=FitH]{hyperref}
\fi
\usepackage{epsfig}
\usepackage{graphicx}

\newtheorem{theorem}{Theorem}

\newtheorem{lemma}{Lemma}

\newtheorem{definition}{Definition}

\newtheorem{claim}{Claim}

\newtheorem{remark}{Remark}
\newenvironment{proof}{\noindent{\bf Proof:~~}}{}
\journal{Theory of Computer Science}

\newcommand{\ie}{\emph{i.e.,\ }}

\newcommand{\namedref}[2]{\hyperref[#2]{#1~\ref*{#2}}}

\newcommand{\sectionref}[1]{\namedref{Section}{#1}}

\newcommand{\theoremref}[1]{\namedref{Theorem}{#1}}

\newcommand{\figureref}[1]{\namedref{Figure}{#1}}

\newcommand{\lemmaref}[1]{\namedref{Lemma}{#1}}

\newcommand{\namedrefeq}[2]{\hyperref[#2]{#1~\mbox{\rm(\ref*{#2})}}}
\newcommand{\equationref}[1]{\namedrefeq{Eq.}{#1}}
\newcommand{\remarkref}[1]{\namedref{Remark}{#1}}
\newcommand{\definitionref}[1]{\namedref{Definition}{#1}}

\newcommand{\hide}[1]{}

\newcommand{\tb}{\makebox[0.4cm]{}}
\newcommand{\due}{\makebox[0.8cm]{}}

\newcommand\A{\mathcal{A}}

\newcommand\C{\mathcal{C}}

\newcommand\N{N}

\newcommand\norm[1]{\left|\left|#1\right|\right|}
\newcommand\normi[1]{\left|\left|#1\right|\right|_\infty}

\newsavebox{\theorembox}
\newsavebox{\lemmabox}
\newsavebox{\conjecturebox}
\newsavebox{\claimbox}
\newsavebox{\factbox}
\newsavebox{\corollarybox}

\newsavebox{\propositionbox}
\newsavebox{\examplebox}
\savebox{\theorembox}{\bf Theorem}

\savebox{\lemmabox}{\bf Lemma}

\savebox{\conjecturebox}{\bf Conjecture}

\savebox{\claimbox}{\bf Claim} \savebox{\factbox}{\bf Fact}

\savebox{\corollarybox}{\bf Corollary}

\savebox{\propositionbox}{\bf Proposition}

\savebox{\examplebox}{\bf Example}

\newcommand{\ignore}[1]{}

\newcommand{\syncAlg}{\textsc{SS-Iterative}\xspace}
\newcommand{\asyncAlg}{\textsc{Async-SS-Iterative}\xspace}
\newcommand{\linenumber}[1]{{\tt #1}}
\newcommand{\lineref}[1]{Line~\linenumber{#1}}
\newcommand{\uu}{\mathbf{u}}
\newcommand{\vv}{\mathbf{v}}
\newcommand{\cc}{\mathbf{c}}
\newcommand{\xx}{\mathbf{x}}

\newcommand{\bb}{\mathbf{b}}
\newcommand{\diag}{\mathop{\mathrm{diag}}}

\def\squarebox#1{\hbox to #1{\hfill\vbox to #1{\vfill}}}

\newcommand{\D}{\mathcal{D}}

\newcommand{\Ir}[1]{I(#1)}
\newcommand{\Or}[1]{O(#1)}
\newcommand{\II}[2]{I_{#1}(#2)}
\newcommand{\OO}[2]{O_{#1}(#2)}

\begin{document}

\begin{frontmatter}



\title{Self-stabilizing Numerical Iterative Computation\footnote{Preliminary version of this work appeared
in the 10th International Symposium on Stabilization, Safety, and Security of Distributed Systems (SSS), Detroit, Nov. 2008.}}


\author{Danny Bickson$^a$, Ezra N. Hoch$^b$, Harel Avissar$^b$ and Danny Dolev$^b$}


\address[IBM]{
IBM Haifa Research Lab\\
Mount Carmel\\Haifa 31905, Israel\\
Email: \{dannybi\}@il.ibm.com\\
}
\address[HUJI]{School of Computer Science and Engineering\\
The Hebrew University of Jerusalem\\
Jerusalem 91904, Israel\\
Email: \{hoch,harela01,dolev\}@cs.huji.ac.il\\
}



\begin{abstract}
Many challenging tasks in sensor networks, including sensor
calibration, ranking of nodes, monitoring, event region detection,
collaborative filtering, collaborative signal processing, {\em etc.}, can
be formulated as a problem of solving a linear system of
equations. Several recent works propose different distributed
algorithms for solving these problems, usually by using linear iterative numerical
methods.

The main problem with previous approaches is that once the problem inputs
change during the process of computation, the computation may output unexpected results.
In real life settings, sensor measurements are subject to
varying environmental conditions and to measurement noise.

We present a simple iterative scheme called \syncAlg for solving
systems of linear equations, and examine its properties in the self-stabilizing
perspective. We analyze the behavior of the proposed scheme under changing input sequences using two different assumptions on the input: a box bound, and a probabilistic distribution.

As a case study, we discuss the sensor calibration problem and
provide simulation results to support the applicability of our
approach.
\end{abstract}


\end{frontmatter}



\section{Introduction}
Many challenging tasks in sensor networks, for example distributed ranking algorithms
of nodes and data items \cite{PPNA08}, collaborative filtering \cite{KorenCF},
localization \cite{SensorLocalization}, collaborative signal
processing \cite{SignalProcessing}, region detection
\cite{RegionDetection}, etc., can be formulated as a problem of solving a linear system of
equations. Several recent works \cite{SensorLocalization},\cite{SignalProcessing},\cite{RegionDetection} propose different distributed
algorithms for solving these problems, usually by linear iterative numerical
methods.

In this work, we extend the settings of the above approaches by
adding another dimension to the problem. Specifically, we are
interested in {\em self-stabilizing} algorithms, that continuously run
and converge to a solution from any initial
state. This aspect of the problem is highly important due to the dynamic
nature of the network and the frequent changes in the
values being tracked.

As a case study, we show that the calibration of local sensors' readings can be formulated as a
linear system of equations $A\xx = \bb$, where $\xx$ represents
the calibrated output reading, $\bb$ represents the local reading, and $A$
represents a weighted communication graph. However, our work is general and
can be applied to any problem that can be formulated as a distributed solution to a
linear system of equations, including the previous works mentioned above.

The main challenge we have faced in this work, is that in the
classical linear algebra literature, $\bb$ is assumed to be constant.
In our settings, the environment is constantly changing and the
computed algorithm never terminates, leading to constantly
changing values of $\bb$. In this paper, we ask the following question: ``Is it possible to
devise a self-stabilizing numerical iterative method?'' We answer
affirmatively, and show that under minor conditions it is possible
to devise a self-stabilizing algorithm that solves a dynamic
system of linear equations, where the input to the system is
constantly changing.

We analyze the behavior of the proposed scheme under changing input sequences
using two different assumptions on the input behavior: a box bound, and a probabilistic bound.
We show that in both cases, once the network stabilizes there is a guarantee on
the quality of the computation.

One of the most efficient distributed approaches for solving a set
of linear equations of the type $A\xx = \bb$ is by using linear
iterative algorithms. Unlike Gaussian elimination, which has a
cost of $O(n^3),$ where $n$ is the number of variables, an
iterative algorithm usually solves a system of linear equations in
time of $O(n^2r)$ where $r$ is the number of iterations, which depends
on the properties of the matrix $A$. In many real life problems, the convergence
is logarithmic in $n$, for example in multiuser detection \cite{Allerton}.
These algorithms are naturally distributed and work well in asynchronous settings. Furthermore,
when converging, the algorithms converge to a solution from {\em
any} initial state. An excellent survey of such methods is found
in \cite{BibDB:BookBertsekasTsitsiklis}.

The main novel contribution of this paper is in analyzing
the self-stabilizing properties of
algorithms from the linear iterative methods domain. For example, in a practical setting, it is
highly unreasonable to assume that sensor readings do not change
over time. However, it is reasonable to assume that at steady
state the change in sensor readings is bounded. Informally, in
this work we show that once the input readings are bounded, the
output solution is bounded as well. This useful observation
enables us to tie together numerical iterative methods and
dynamically changing environments in a self-stabilizing manner.

To the best of our knowledge, this is the
first work tackling this challenging problem. We believe that our
approach can have numerous applications in the field of
distributed self-stabilizing computation.

Other works discuss fault tolerance aspects of distributed
computation. For example, overcoming faults in sensors by
averaging the input was investigated in~\cite{SensorFailures},
using a centralized algorithm. Quantifying faulty nodes'
effect on the system's output is discussed in~\cite{BinaryRept}
and~\cite{MultiValRept}. These papers consider bounded input paths
and their effect on the stability of the output. In
\cite{RepCons-PODC08} infinite input paths are considered under
the assumption that only specific sensors' input may change. All
three papers consider discrete input values, as opposed to a
continuous set of input values discussed in this paper.

The paper is constructed as follows.
\sectionref{sec:model} defines the model and presents the problem definition. \sectionref{sec:Errors} discusses measurement error types and characterization. \sectionref{sec:solution}
presents our novel algorithm \syncAlg. \sectionref{sec:analysis}
analyzes our algorithm and gives bounds on the convergence rate, assuming a box bound on the inputs. \sectionref{sec:async} extends our construction to the asynchronous case. \sectionref{sec:Probabilistic} extends the
analysis of \syncAlg under a probabilistic model when the input is characterized using a normal distribution.
\sectionref{sec:simlation} presents experimental results of
running \syncAlg using sample topologies. We conclude in
\sectionref{sec:discussion}.

\section{Model and Problem Definition}\label{sec:model}
Consider a distributed system of sensors measuring real-world
data. Sensors are located in different areas; for example, the
senors are spread throughout a building and they measure the
temperature to adjust the heating or cooling. We would like the
collected data to be as reliable as possible, reflecting closely the changing
environmental conditions. One of the obstacles we face when designing
algorithms that collect data from a sensor network are measurement errors.
There are two main types of inaccuracies of sensors' measurements: noisy environment
and sensing equipment which is not calibrated.
In our setting sensors are allowed to
communicate among themselves
and
to use data from other nodes for
calibrating their
environmental readings.
Furthermore, we would like our
calibration algorithm to have fault-tolerance properties.
Specifically, we are interested in self-stabilizing algorithms
\cite{DolevSSBook}, which converge to an optimal solution from
any initial state. Observe that self-stabilization helps also in deploying the sensors.
There is no need to explicitly synchronize the sensors, once enough of them
are deployed and begin functioning the results will converge to the expected value.

We model the sensor calibration problem as follows. Given a directed
communication graph $G =(V,E)$, $V$ is the set of nodes $V =
\{p_1, \dots, p_n\}$, $E$ is the set of weighted edges connecting
them (weights can be negative). Edge weights are used to model the directional dependence between nodes' outputs; \ie if $w_{p_i,p_j}=0$ then there is no edge between $p_i$ and $p_j$ and their output is not directly dependent on each other. In addition, it is possible to have a non-zero self
connected edge, $w_{p_i,p_i} \neq 0$, which represents the weight of $p_i$'s own input. Note that the graph
does not have to be symmetric, which means that it is possible that $w_{p_i,p_j} \ne w_{p_j,p_i}$.
We denote $\N(p_i)$ as the set of node $p_i$'s neighbors, excluding $p_i$. 

Initially, we assume a synchronous system: during a single round of communication, any pair of connected
nodes may send a single message on each directed edge. In each round
$r$, each node $p_i$ has a scalar input value $\II{p_i}{r}$, which
represents the local reading of the sensor.\footnote{For simplicity
of notations we use scalar variables in the paper. An extension to
the vector case (where each sensor measures a set of measurements) is immediate.} In addition, $p_i$ outputs its
output value, which is denoted by $\OO{p_i}{r}$; both inputs and outputs are from the domain of real numbers. Denote by $\Ir{r}$ the
input vector of the entire system at round $r$, and by $\Or{r}$
the output vector of the system at the end of round $r$.
In \sectionref{sec:async} we relax the assumption of synchronous rounds
and provide a variant of the algorithm that works in asynchronous settings.

The schematic operation of each node $p_i$ at round $r$ is
composed of the following steps: (a) read the value of
$\II{p_i}{r}$; (b) send messages; (c) receive messages; and (d) do
some processing and output $\OO{p_i}{r}$. Then a new round is
started, and the nodes continue so forever.

\begin{definition}
  A {\bf configuration} $\C$ of the system at round $r$ consists of the state
  of each node prior to performing any operation at round $r$;
  this configuration is denoted by $\C(r)$.
\end{definition}

\begin{definition}
  An {\bf input sequence} $\mathcal{I}$ of length $\ell$ is a list of  $\ell$ vectors such that each $\vv \in \mathcal{I}$
  is a possible input vector of the system (\ie $\vv \in \mathfrak{D}$, the domain of allowed values).
  An {\bf output sequence} $\mathcal{O}$ of length $\ell$ is a list
  such that each $\vv \in \mathcal{O}$
  can potentially be an output vector of the system.
\end{definition}

\begin{definition}
  A step from configuration $\C$ to configuration $\C'$ on input vector
  $\vv$ is {\bf legal} if $\C'$ is reached from $\C$ by the
  system when having $\vv$ as the input vector.
  $\uu$ is {\bf produced} by a legal step if
  $\uu$ is the output vector of the system resulting from such a legal step.
\end{definition}

\begin{definition}
  A {\bf run} of a system on input sequence $\mathcal{I} = \{\vv(1), \dots, \vv(\ell)\}$ starting
  from configuration $\C(r)$ is the sequence $\C(r), \Or{r}, \C(r+1), \Or{r+1}, \dots$
  s.t. for any $i \geq 0$: the step from $\C(r+i)$ to $\C(r+i+1)$ on input vector
  $\vv(i+1)$ is legal, and  $\Or{r+i}$ is produced by that
  legal step. The system is said to produce the output sequence
  $\mathcal{O} = \{\Or{r}, \dots, \Or{r+\ell-1}\}$.
\end{definition}

In the special case when the sensor observations (the input to the
system) are fixed, the output decision of the sensors should converge to
a solution that preserves the linear relations among node inputs
and outputs. More formally, consider an input sequence
$\mathcal{I}$ of identical input vectors; \ie $\mathcal{I} =
\{\vv, \vv, \vv, \dots\}$. It is desired that for such an input a run
from any configuration $\C$ on $\mathcal{I}$ would end up
producing an output sequence $\mathcal{O} = \{\uu(1), \uu(2),
\dots\}$ such that $\norm{\uu(i)-\uu} \rightarrow 0$ as $i
\rightarrow \infty$, for a $\uu$ that solves the following linear
system of equations:
\begin{equation}\label{eq:Opdef}
  \uu_i =
  w_{p_i,p_i} \vv_i + \sum_{p_j \in \N(p_i)}w_{p_i, p_j}
  \uu_j\;.
\end{equation}
We assume that the above set of equations are uniquely solvable, denoting
$\uu$ as the solution to $\vv$.

The following definition bounds the change in sensor observations:
\begin{definition}\label{def:concentrated}
  An input sequence $\mathcal{I} = \{\vv(1), \vv(2), \dots, \vv(\ell)\}$ is $\delta$-{\bf bounded} around
  $\vv$ if for every $i$, $1 \leq i \leq \ell$, it holds that $\norm{\vv(i)-\vv}_\infty \leq
  \delta$.\footnote{$\norm{\xx}_\infty = \max_i \{ |\xx_i| \}$.}
\end{definition}

\definitionref{def:concentrated} states that a sequence $\mathcal{I}$ is $\delta$-bounded if all
the vectors in $\mathcal{I}$
are bounded within an $n$ dimensional hypercube with an edge $2
\delta$, centered around a point $\vv$. We note that once changes
in the input are not bounded, then no efficient algorithm
(especially in a network that is sparsely connected) can calculate
the output fast enough. Thus it is easy to construct an example where the diameter of the
communication graph is $\D$, for some system of equations it would
take at least $\D$ rounds for the information exchange for input readings at
one side of the network to propagate to the other side of the
network.

\begin{definition}\label{def:runconverges}
  Let $\mathcal{I}$ be an input sequence that is $\delta$-bounded
  around $\vv$ and let $\uu$ be the solution to input $\vv$.
  A run from configuration $\C$ on input sequence $\mathcal{I}$ $\epsilon$-{\bf converges} to its solution if
  the produced output sequence $\mathcal{O} = \{\uu(1), \uu(2), \dots\, \uu(\Delta t)\}$ satisfies that
  $\normi{\uu(\Delta t)-\uu)}\leq \epsilon(\Delta t, \delta, \C)$; where $\epsilon$
  is a function of $\Delta t, \delta$ and $\C$.
\end{definition}

\definitionref{def:runconverges} requires that if - starting from
configuration $\C$ - the inputs are in an $n$ dimensional hypercube of
radius $\delta$ around $\vv$ then the output at time $\Delta t$ is
bounded within some $n$ dimensional hypercube around $\uu$ with
radius $\epsilon(\Delta t, \delta, \C)$. We aim at an
$\epsilon(\Delta t, \delta, \C)$ that decreases as $\Delta t$
increases, as long as the inputs are bounded by the same
$\vv$-centered, $\delta$-radius hypercube. Clearly, for $\delta
>0$, $\epsilon(\Delta t, \delta, \C) > 0$ for any $\Delta t$. That is, there
is some minimal radius $\delta' >0$ around $\uu$ and we cannot
ensure a tighter bound.

The above definition considers a single initial configuration, and
a single input sequence $\mathcal{I}$. We are interested in an
algorithm that works for all initial configurations and all input
sequences.

\begin{definition}\label{def:alwaysconv}
  An algorithm $\A$ $\epsilon$-{\bf converges} for $\delta$-bounded input sequence $\mathcal{I}$ if every run (from any configuration) on
  $\mathcal{I}$,
  $\epsilon$-converges to its solution.  An algorithm $\A$ $\epsilon$-{\bf always converges} if for
  every $\delta$-bounded input sequence $\mathcal{I}$, $\A$ $\epsilon$-converges.
\end{definition}

\definitionref{def:alwaysconv} formally defines the problem at
hand, as an algorithm $\A$ that {\em always converges} has the
desired self-stabilizing property: for any system state, once the
sensors' readings changes are bounded, the change in output of the
entire system is bounded as well.

Our goal is to find an algorithm $\A$ that is $\epsilon$-always
converging for a provably ``good'' $\epsilon$. Moreover, we aim at
having $\A$ efficient also in its message complexity and
the simplicity of the code, allowing lightweight sensors to actually
implement it.

\section{Measurement Errors}\label{sec:Errors}
As mentioned in the introduction, sensor readings suffer from two types of inaccuracies: measurement errors
and calibration errors. In this section we elaborate on two common types of measurement errors.
This introduction serves as a motivation for bounding the errors using a normal distribution,
a construction that is discussed in detail in \sectionref{sec:Probabilistic}.

Most physical sensors tend to suffer from noise, which is defined as fluctuations in external factors to the stream of target information (signal). Most sensor noise can be attributed to a small number of common noise patterns, this noise is usually due to the nature of physical substance being measured and the nature of the measurement equipment. Two of the most common noise patterns are Shot Noise and Nyquist Noise as defined at \cite{RFCD97}.

\begin{definition}\label{def:ShotNoise}
{\bf Shot Noise} is a type of noise common to the particle-like nature of the charge carriers. It is often thought that a DC current flow in any semiconductor material is constant at every instant. In fact, however, since current flow is made up of individual electrons and holes, it is only the time average flow of these charge carriers that is seen as a constant current. Any fluctuation in the number of charge carriers at any instance produces a random current change in the instance.
\end{definition}
 As these particles (usually photons or electrons) distribute with Poisson distribution the noise itself also follows this distribution. Poisson distribution is nearly Gaussian and can easily be approximated by a Normal distribution.

\begin{definition}\label{def:NyquistNoise}
{\bf Nyquist Noise} (also Thermal noise): In any conducting medium whose temperature is above absolute zero (0 Kelvin), the random motion of charge carriers within the conductor produces random voltages and currents. These voltages and currents produce noise.
\end{definition}
 This noise is nearly white by nature (equally distributed through the spectrum). It is characterized by a normal distribution of its amplitude.

Both noise types discussed above are nearly Gaussian \cite{GIE}. Therefore tolerance to Gaussian noise is an important issue when examining the sensor calibration problem, specifically, and for systems in general. Furthermore, the central limit theorem states that in the average case any i.i.d sample is bound to be distributed as a Gaussian surrounding the mean. This makes Gaussian tolerance a powerful tool when dealing with dynamic input sequences.

\section{Our Proposed Solution}\label{sec:solution}
%
%
Getting back to the sensor calibration example, let $W$ be the matrix that has $w_{i,j}$ as entries,
\equationref{eq:Opdef}. It be written in linear algebra
notation, (s.t. it applies to all nodes simultaneously):
\begin{equation} \label{eq:matrixform}
  W\uu = \vv\;.
\end{equation}

If we consider a non-self-stabilizing system in which the inputs
do not change (that is, the input is fixed, as some $\vv$), then
\equationref{eq:matrixform} can be seen as $A\xx=\bb$, where $A$
and $\bb$ are given. In such a case, we are interested in finding
the value of $\xx$, which is a vector of $n$ unknown variables.
However, we are interested in the case where $\vv$ changes over
time, and thus \equationref{eq:matrixform} does not describe the
problem properly, but rather helps in understanding the motivation
for our solution.

Rearranging \equationref{eq:Opdef} we get:
\begin{equation}\label{eq:updaterule}
  \OO{p_i}{r+1} =
  w_{p_i,p_i} \II{p_i}{r+1} + \sum_{j \neq i}w_{p_i, p_j}
   \OO{p_j}{r}\;.
\end{equation}

Clearly, for the case of $\delta=0$, a $0$-bounded input sequence
$\mathcal{I}$, if $(\OO{p_i}{r+1}-\OO{p_i}{r}) \longrightarrow 0$
as $r \rightarrow \infty$ then \equationref{eq:updaterule}
converges to the solution of \equationref{eq:Opdef}. Thus, if the
update rule of \equationref{eq:updaterule} is executed
simultaneously by all nodes, and for all of the nodes
$(\OO{p_i}{r+1}-\OO{p_i}{r}) \longrightarrow 0$, then it also
solves \equationref{eq:matrixform}. That is, if each node locally
executes \equationref{eq:updaterule} then the global solution is
reached. This observation motivates algorithm \syncAlg in
\figureref{figure:syncAlg}.

\begin{figure*}[h]
\begin{center}
\begin{minipage}{4.8in}
\hrule \hrule \vspace{1.7mm} \footnotesize
\setlength{\baselineskip}{3.9mm} \noindent Algorithm \syncAlg
 \vspace{1mm} \hrule \hrule
\vspace{1mm}

\linenumber{01:} Each round {\bf do}:\hfill\textit{/* executed on
node $p_i$
*/}\\

\makebox[0.93cm]{} \textit{/* send current value of $O_{p_i}$ to all neighbors */}\\
\linenumber{02:} \tb {\bf for} each $p_j \in \N(p_i)$  \\
\linenumber{03:} \due {\bf send} $O_{p_i}$ to $p_j$; \\

\makebox[0.93cm]{} \textit{/* update $O_{p_i}$ according to values sent by neighbors */}\\
\linenumber{04:} \tb {\bf set} $O_{p_i} := w_{p_i,p_i} I_{p_i};$  \\
\linenumber{05:} \tb {\bf for} each value $O_{p_j}$ received:  \\
\linenumber{06:} \due {\bf update} $O_{p_i} := O_{p_i}  + w_{p_i,p_j} O_{p_j}$;\\

\linenumber{07:} {\bf od}.

\normalsize \vspace{1mm} \hrule\hrule
\end{minipage}
 \caption{A self-stabilizing iterative algorithm.}\label{figure:syncAlg}
 \end{center}
\end{figure*}

\begin{remark}
  In \syncAlg there is no notion of the ``current round number
  $r$''. That is, $p_i$ reads and writes to the variables $I_{p_i}$ and $O_{p_i}$
  without being ``aware'' of $r$. When we discuss the algorithm ``from the outside'', we will consider $\II{p_i}{r}$ and $\OO{p_i}{r}$ instead of
  just $I_{p_i}, O_{p_i}$.
\end{remark}

Consider $p_i$ as running at some round $r+1$. When $p_i$ performs
\lineref{03}, it sends the value of $O_{p_i}$. The last time
$O_{p_i}$ was updated was at \lineref{04} and \lineref{06} of
round $r$. Thus, the value sent by $p_i$ at round $r+1$ is
actually $\OO{p_i}{r}$. Therefore, the values received from $p_j$
by $p_i$ and used to update $\OO{p_i}{r+1}$ are $\OO{p_j}{r}$.
However, the value read by $p_i$ in \lineref{04} is the value of
$\II{p_i}{r+1}$. Concluding that $p_i$ updates $\OO{p_i}{r+1}$
exactly according to \equationref{eq:updaterule}.

\begin{remark}\label{rem:partofcode}
  Each node $p_i$ must know the values of $w_{p_i, p_j}$ as ``part of the
  code''. Thus, these values cannot be subject to transient
  faults.
\end{remark}

\section{Analysis of \syncAlg}\label{sec:analysis}
Paper \cite{BibDB:BookBertsekasTsitsiklis} shows that the update rule of
\equationref{eq:updaterule} can be written in linear algebra form
as
\begin{equation}\label{eq:matrixUpdate}
  \Or{r+1} = A\Ir{r+1} + B\Or{r}\;,
\end{equation}
where $A$ is a diagonal matrix with $w_{p_i,p_i}$ in the main
diagonal, and $B_{ij} = w_{p_i, p_j}$ for $i\neq j$
\begin{align}
 A &\triangleq (\diag(W))^{-1}\,, &  B &\triangleq (I_{n \times n} - \diag(W)^{-1}W) \label{AB}\,,
\end{align}
 where $I_{n \times n}$ is the identity matrix. Using this update rule to solve a set of linear equations iteratively is  known as the Jacobi
algorithm.


Let $\mathcal{I}$ be an input sequence of length $\ell$ that is
$\delta$-bounded around vector $\vv$. At some round $r,$  $\mathcal{I} =
\Ir{r}, \Ir{r+1}, \dots, \Ir{r+\ell-1}$. Note
that \syncAlg saves a single scalar variable at each node, and thus the
configuration of round $r+1$ can be defined by the value of
$\Or{r}$ at round $r$. Consider \syncAlg's run, starting from an
arbitrary configuration at round $r$. We aim at showing that
$\Or{r+\Delta t}$ is bounded by a hypercube centered at $\uu$. Denote
by $\cc(\Delta t) \triangleq \Or{r+ \Delta t}-\uu$. If we show
that $\normi{\cc(\Delta t)}$ is bounded (as $\Delta t$ increases),
then $\Or{r+\Delta t}$ is within a bounded hypercube centered at $\uu$.
Consider $\cc(1)$:
\begin{eqnarray} \label{eq:baseInd}
  \cc(1) & = & \Or{r+1}-\uu \nonumber\\
         & = & A\Ir{r+1} + B\Or{r} -(A\vv+B\uu) \nonumber\\
         & = & A(\Ir{r+1}-\vv) + B(\Or{r}-\uu) \nonumber\\
         & = & A(\Ir{r+1}-\vv) + B\cc(0)\;.
\end{eqnarray}
Since $\mathcal{I}$ is a $\delta$-bounded input sequence around
$\vv$, each $\Ir{r+\Delta t}$ can be denoted as
$\vv+\D(r+\Delta t)$ s.t. $\D(r+\Delta t) \in
\mathbb{R}^n$ is a vector, and $\norm{\D(r+\Delta t)}_\infty
\leq \delta$. 


\begin{claim}
At round $r+\Delta t$, it holds that \[ \cc(\Delta t)= \sum_{j=0}^{\Delta t-1} B^jA\D(r+\Delta t-j) + B^{\Delta t}\cc(0). \]
\end{claim}
\begin{proof}
  Proof by induction. The base of the induction was shown for $\cc(1)$; see \equationref{eq:baseInd}. Assume that the claim holds for $\Delta t =
  k$. Thus, $\cc(k)= \sum_{j=0}^{k-1} B^jA\D(r+k-j) + B^k\cc(0)$. By the update rule in
  \equationref{eq:matrixUpdate}, we have that $\Or{r+k+1} = A\Ir{r+k+1}
  + B\Or{r+k}$. Combining the two equations implies
  \begin{eqnarray*}\label{eq:norm}
    \cc(k+1) & = & \Or{r+k+1} - \uu \\
           & = & A\Ir{r+k+1} + B\Or{r+k} - (A\vv+B\uu)\\
           & = & A(\Ir{r+k+1}-\vv) + B(\Or{r+k} - \uu)\\
           & = & A\D(r+k+1) + B\cc(k)\\
           & = & A\D(r+k+1) + \sum_{j=0}^{k-1} B^{j+1}A\D(r+k-j) + B^{k+1}\cc(0)\\
           & = & A\D(r+k+1) + \sum_{j=1}^{k} B^jA\D(r+k+1-j) + B^{k+1}\cc(0)\\
           & = & \sum_{j=0}^{k} B^jA\D(r+k+1-j) + B^{k+1}\cc(0)\;.
  \end{eqnarray*}
  Thus, if the claim holds for $\Delta t=k$ it also holds for $\Delta t=k+1$;
  and we have that the claim holds for all $\Delta t \geq 0$.
  \qed
\end{proof}

\begin{definition}
  A matrix $M_{n \times n}$ is {\bf diagonally dominant} if
  $|M_{ii}| > \Sigma_{j \ne i}^n |M_{ij}| $. A matrix $M_{n \times n}$ is {\bf normalized} diagonally dominant (normalized, for short) if
  $M$ is diagonally dominant, and $|M_{ii}| \geq 1$.
\end{definition}

\begin{lemma}\label{claim:diag}
  For a normalized diagonally dominant matrix $W$, it holds that
  $\normi{A} \leq 1$ and $\normi{B} < 1,$ where $A,B$ are defined in \equationref{AB} and $\normi{A} \triangleq \max_{\xx \neq 0}\frac{\normi{Ax}}{\normi{x}}$.
\end{lemma}
\begin{proof}
  $A$ is zero except for its
  main diagonal for which $A_{i,i} = w_{p_i, p_i} =\frac{1}{w_{i,i}}$. Since $|W_{i,i}|
  \geq 1$, we have that $|A_{i,i}| \leq 1$. Thus, it holds that $\norm{A\xx}_\infty \leq \norm{\xx}_\infty$.
  Furthermore,
  $\max_{\xx \neq 0}\frac{\norm{A\xx}_\infty}{\norm{\xx}_\infty}
  \leq
  1$, \ie $\norm{A}_\infty \leq 1$. Regarding $B$, $B_{i,j} = w_{p_i, p_j}$ for $i
  \neq j$ and 0 for $i=j$. Since $W$ is assumed to be normalized diagonally dominant, we have
  that $\sum_{j \ne i} |W_{i,j}|< |W_{i,i}|$, thus $\sum_{j \ne i}
  |w_{p_i,p_j}|< 1$. Therefore, $\sum_j |B_{i,j}| = \sum_{j \ne i} |w_{p_i,p_j}|<
  1$ for all $i$. In total, for any $\xx$ we have
  $\norm{B\xx}_\infty < \norm{\xx}_\infty$, leading to $\norm{B}_\infty <
  1$.
\qed\end{proof}

If $W$ is a diagonally dominant matrix then node $p_i$'s own input effects $p_i$'s output more than the sum of all of $p_i$'s neighbors outputs. That is, the weight of $p_i$'s input is at least the sum of weights of $p_i$'s neighbors outputs.

\begin{theorem}\label{theorem:main}
  Given a normalized diagonally dominant and invertible $W$, there are constants $c_1,c_2,$ where $c_1 > 0,$ and $ 1 > c_2 > 0$,
  such that \syncAlg $\epsilon$-always converges with
  $\epsilon(\Delta t, \delta, \C) = \delta c_1 + (c_2)^{\Delta t} \normi{\Or{r}-\uu}$.
\end{theorem}
\begin{proof}
  By \lemmaref{claim:diag} it holds that $\norm{A}_\infty \leq 1$ and $\norm{B}_\infty <
  1$. Consider a $\delta$-bounded input sequence $\mathcal{I}$ around $\vv$, and
  \syncAlg's run starting from an arbitrary state $\Or{r}$.
  We are interested in the behavior of $\normi{\cc(\Delta t)}$:
\begin{eqnarray}\label{eq:normck}
  \normi{\cc(\Delta t)} & = & \normi{\sum_{j=0}^{\Delta t-1} B^jA\D(r+\Delta t-j) + B^{\Delta t}\cc(0)} \nonumber\\
              & \leq & \normi{\sum_{j=0}^{\Delta t-1} B^jA\D(r+\Delta t-j)} + \normi{B^{\Delta t}\cc(0)} \nonumber\\
              & \leq & \sum_{j=0}^{\Delta t-1} \normi{B}^j\normi{A\D(r+\Delta t-j)} + \normi{B}^{\Delta t}\normi{\cc(0)} \nonumber\\
              & \leq & \delta \normi{A}\sum_{j=0}^{\Delta t-1} \normi{B}^j + \normi{B}^{\Delta t}\normi{\cc(0)} \nonumber \\
              & = & \delta \normi{A} \frac{1-\normi{B}^{\Delta t}}{1-\normi{B}} + \normi{B}^{\Delta t}\normi{\cc(0)}\;.
\end{eqnarray}
For an input sequence $\mathcal{I}$ that is $\delta$-bounded around
$\vv$, denote by $\uu$ the solution to the original system of
equations $W\uu = \vv$. By \equationref{eq:norm},
\begin{equation*}
   \normi{\Or{r+\Delta t}-\uu} \leq \delta \normi{A} \frac{1-\normi{B}^{\Delta t}}{1-\normi{B}} +
  \normi{B}^{\Delta t}\normi{\cc(0)}.
\end{equation*}
  Since $\normi{B} < 1$, we have that
  $\frac{1-\normi{B}^{\Delta t}}{1-\normi{B}}\leq
  \frac{1}{1-\normi{B}}$ and by setting $c_1 =
  \frac{\normi{A}}{1-\normi{B}}$ it holds that $\normi{A}
  \frac{1-\normi{B}^{\Delta t}}{1-\normi{B}} \leq c_1$. By setting $c_2 =
  \normi{B}$ and recalling that $\cc(0) = \Or{r}-\uu$ we are done.
\qed\end{proof}

\theoremref{theorem:main} states sufficient conditions s.t. \syncAlg $\epsilon$-always
converges. Moreover, the algorithm \syncAlg is lightweight, as it
requires nodes to send only a single value to every neighbor on each
round.

\section{Extension to the Asynchronous Model}\label{sec:async}
Our second novel contribution is in extending our model to support
asynchronous communications. In a large sensor network, it is
unreasonable to assume that the sensors operate in synchronous
rounds. Furthermore, as known from the linear iterative algorithms
literature, algorithms
sometimes
converge in less asynchronous rounds (when compared to synchronous rounds).

When considering the asynchronous model, it is more convenient to
discuss shared-memory as means of communication.\footnote{In
\cite{DolevSSBook} it is shown how to convert an algorithm based
on shared-memory to a message-passing algorithm with links of
bounded capacity.} Thus, assume that for each directed edge
between $p_i, p_j$ there is a read-write register $R_{p_i,p_j}$
that is written by $p_i$ and read by $p_j$.

An asynchronous run is an infinite sequence of configurations $\C_0
\rightarrow \C_1 \rightarrow \dots$ such that some process $p$
performs an atomic step between configuration $\C_i$ and
$\C_{i+1}$. An atomic step consists of reading or writing from a
single register. Notice that in the current model a configuration
consists of all the registers and of the local variables at the
different nodes.

In this section we again prove that starting from an arbitrary
configuration, when the changes to the inputs are bounded, the outputs are
bounded as well. We consider each configuration $\C_r$ to be
assigned a vector input $\Ir{r}$ such that if node $p_i$ reads the
input when performing an atomic step on $\C_r$ it reads the value
of $\II{p_i}{r}$. Equivalently, the output vector of configuration
$\C_r$ is $\Or{r}$.

\figureref{figure:asyncAlg} outlines \asyncAlg, which is a direct
translation of \syncAlg to the shared-memory model.

\begin{figure*}[h]
\begin{center}
\begin{minipage}{4.8in}
\hrule \hrule \vspace{1.7mm} \footnotesize
\setlength{\baselineskip}{3.9mm} \noindent Algorithm \asyncAlg
 \vspace{1mm} \hrule \hrule
\vspace{1mm}

\linenumber{01:} Forever {\bf do}:\hfill\textit{/* executed on
node $p_i$
*/}\\

\makebox[0.93cm]{} \textit{/* write current value of $O_{p_i}$ to all neighbors */}\\
\linenumber{02:} \tb {\bf for} each $p_j \in \N(p_i)$:  \\
\linenumber{03:} \due {\bf write} $O_{p_i}$ to $R_{p_i,p_j}$; \\

\makebox[0.93cm]{} \textit{/* update $O_{p_i}$ according to values of neighbors */}\\
\linenumber{04:} \tb {\bf set} $O_{p_i} := w_{p_i,p_i} I_{p_i};$  \\
\linenumber{05:} \tb {\bf for} each $p_j \in \N(p_i)$:  \\
\linenumber{06:} \due {\bf read} $R_{p_j,p_i}$ into $temp$;  \\
\linenumber{07:} \due {\bf update} $O_{p_i} := O_{p_i}  + w_{p_i,p_j} temp$;\\

\linenumber{08:} {\bf od}.

\normalsize \vspace{1mm} \hrule\hrule
\end{minipage}
 \caption{A self-stabilizing iterative algorithm for asynchronous networks.}\label{figure:asyncAlg}
 \end{center}
\end{figure*}

\asyncAlg consists of two phases: in the first, the previous value of $O_{p_i}$ is
written to all its neighbors' registers. In the second phase $p_i$ calculates
its new value of $O_{p_i}$ by reading the registers of all its
neighbors.

We consider only ``fair'' runs, in which each node performs an
atomic step infinitely many times. Thus, each node performs both
phases infinitely many times. A round is defined to be the shortest
prefix of a run such that each node has performed at least once lines \linenumber{02}-\linenumber{07} in the
algorithm. We number each atomic step and each round. Note that a round consists of many
atomic steps.

We model a fair run as follows. Each node $p_i$ performs
infinitely many atomic steps, and participates in infinitely many
rounds. Notice that the registers $p_i$ reads in round $k+1$ have
all been last written to, no earlier than during round $k$. Since a
round consists of each node performing all the steps in the
algorithm, each node $p_i$ manages to read all of its neighboring
registers and to write to all of its neighboring registers every
round. Thus, there is some atomic step $r$ (during round $k+1$)
such that:
%
%
\begin{equation*}
  \OO{p_i}{r} =
  w_{p_i,p_i} \II{p_i}{r'} + \sum_{j \neq i}w_{p_i, p_j}
   \OO{p_j}{r_j'}\;,
\end{equation*}
where $r', r_j'$ (for all $p_j \neq p_i$) are smaller
than $r$ and are from at least round $k$.

Let $\uu$ be such that $\uu = A\vv + B\uu$, and let the inputs be
from a $\delta$-bounded input sequence around $\vv$. Denote $\cc(r) = \Or{r} - \uu$ and $z = \max_i |\cc_{p_i}(0)|$.

\begin{theorem}\label{theorem:main2}
  Given a normalized diagonally dominant and invertible $W$, and while considering only fair runs,
  there are constants $c_1,c_2,$ where $c_1 > 0,$ and $ 1 > c_2 > 0$,
  such that \asyncAlg $\epsilon$-always converges with
  $\epsilon(\Delta t, \delta, \C) = \delta  c_1 + (c_2)^{\Delta t}  z$; where $\Delta t$ counts the asynchronous rounds of a fair run.
\end{theorem}

\begin{proof}
 Notice that if $p_i$ did not perform the
$r$th atomic step then $\OO{p_i}{r} = \OO{p_i}{r-1}$ and therefore
$\cc_{p_i}(r) = \cc_{p_i}(r-1)$. Consider the value of
$\cc_{p_i}(r)$ when $p_i$ did perform the $r$th atomic step
(during round $k+1$).
\begin{eqnarray*}
  \cc_{p_i}(r) & = & \OO{p_i}{r} - \uu_{p_i}\\
             & = & w_{p_i,p_i} \II{p_i}{r'} + \sum_{j \neq i}w_{p_i, p_j}  \OO{p_j}{r_j'} - w_{p_i,p_i}\vv_{p_i} - \sum_{j \neq i}w_{p_i, p_j}\uu_{p_i}\\
             & = & w_{p_i,p_i} (\II{p_i}{r'}-\vv_{p_i}) + \sum_{j \neq i}w_{p_i, p_j}  (\OO{p_j}{r_j'}-\uu_{p_i}) \\
             & = & w_{p_i,p_i} (\II{p_i}{r'}-\vv_{p_i}) + \sum_{j \neq i}w_{p_i, p_j}  \cc_{p_j}(r_j')\;,
\end{eqnarray*}
where $r'$ and the different $r_j'$ are smaller than $r$ and are
all from round $k$ or round $k+1$.

By using \lemmaref{claim:diag} we get:
\begin{eqnarray*}
  |\cc_{p_i}(r)| & \leq & |w_{p_i,p_i} (\II{p_i}{r'}-\vv_{p_i})| + \max_{p_j}|c_{p_j}(r_j')| \sum_{j \neq i}|w_{p_i, p_j}| \\
               & \leq & |w_{p_i,p_i} (\II{p_i}{r'}-\vv_{p_i})| + \normi{B}|\cc_{p_{max}}(r_{max})| \\
               & \leq & \delta + \normi{B}|\cc_{p_{max}}(r_{max})|\;,
\end{eqnarray*}
for some $p_{max}$ and $r_{max} \leq r$ that is from round $k$ or
$k+1$.

Therefore, for any $p_i$ during round $k+1$ there is a list of
length $\ell \geq k$ of nodes $p_1, p_2, \dots, p_\ell$ and a
sequence of length $\ell$ of atomic steps $r_1 > r_2 > \dots >
r_{\ell} = 0$, such that
\begin{eqnarray*}
  |\cc_{p_i}(r)| & \leq & \delta + \normi{B}|\cc_{p_1}(r_1)|\\
                     & \leq & \delta + \normi{B}(\delta + \normi{B}|\cc_{p_2}(r_2)|)\\
                     & = & \delta  (1 + \normi{B}) + \normi{B}^2|\cc_{p_2}(r_2)|\\
                     & \leq & \delta  \sum_{z=0}^{\ell-1}\normi{B}^z + \normi{B}^\ell|\cc_{p_\ell}(r_\ell)|\\
                     & = & \delta  \frac{1-\normi{B}^\ell}{1-\normi{B}} + \normi{B}^\ell|\cc_{p_\ell}(0)|\;.
\end{eqnarray*}

Denote by $c_1 \triangleq
\frac{1}{1-\normi{B}}$, and $c_2 \triangleq \normi{B}$. We have
that for node $p_i$ performing the $r$th atomic step during round
$k$ it holds that $|\cc_{p_i}(r)| \leq \delta  c_1 + c_2^\ell
 z \leq \delta  c_1 + c_2^k  z$.
\qed\end{proof}

In fair runs,
there are infinitely many rounds $k$, thus, as $l$ and $r$ go to
infinity, we have that $\normi{\Or{r}}$ is bounded by a hypercube of
length $\delta  c_1$ around $\uu$.

\section{Probabilistic Bound}\label{sec:Probabilistic}
In this section we analyze the behavior of the \syncAlg algorithm, when the input readings are normally distributed
due to noisy readings. 

Let $\mathcal{I}$ be an input sequence of length $\ell$ that is normally distributed around a mean vector $\vv$. That is, $\mathcal{I} =\Ir{r}, \Ir{r+1}, \dots, \Ir{r+\ell-1}$, where each is sampled i.i.d from a Gaussian distribution with mean $\vv$ and covariance $\Sigma_{\vv}$. Consider \syncAlg's run, starting from an arbitrary configuration at round $r$. We aim at showing that $\Or{r}$ is equivalent to a sequence $\mathcal{O}$ sampled i.i.d from a Gaussian distribution $\mathcal{O} \sim \mathcal{N}(\uu, \Sigma_{\uu})$ where $\uu$ is the solution to $\vv$.

Since we know the distribution $\mathcal{I}$, we can calculate $\Or{r}$ using \equationref{eq:matrixUpdate} after $\Ir{r}$ is sampled. $\Or{r}$, our sequence's starting point, can take arbitrary values as we will show that \syncAlg  converges from any starting position. 

We can now deduce that after $k$ rounds of sampling the input distribution the output $\Or{k+1}$ is distributed as:
\begin{equation}\label{eq:perround}
  \Or{k+1} =
  A\Ir{k+1}+B \Or{k} \sim \mathcal{N}(A \vv+B \Or{k}, A \Sigma_{\vv} A^T)\,.
\end{equation}
This formulation is, however, dependant upon $\Or{k}$, or upon the specific samples of $\Ir{\cdot}$. The output of each round is normally distributed with mean $A \vv+B \Or{k}$ and with covariance matrix $A \Sigma_{\vv}A^T$.


To find $\mathcal{O}$ we examine the recursive formula for $\Or{k+1}$. For simplifying notations, we mark the starting
round $r$ as round $0$.
\begin{eqnarray*}
  \Or{k+1} & = & A \Ir{k+1}+B \Or{k}\\
                     & = & A \Ir{k+1}+B A \Ir{k}+B^2\Or{k-1}\\
                     & \vdots & \\
                     & = & A \Ir{k+1}+B A\Ir{k}+\ldots+B^k (A \Ir{1}+B\Or{0})\,.
\end{eqnarray*}
Now we take $k\rightarrow\infty$ to find the distribution $\mathcal{O}$:
\begin{eqnarray}\label{eq:recursiveForm}
  \mathcal{O} \propto \lim_{k\rightarrow\infty}(I_{n \times n}+B+B^2+\ldots+B^{k}) A \mathcal{I}+\lim_{k\rightarrow\infty}B^{k+1} \Or{0}\,.
\end{eqnarray}
By \lemmaref{claim:diag} it holds that $\norm{B}_\infty < 1, \lim_{k\rightarrow\infty}B^{k} = 0$. Using the Maclaurin series we get:
\begin{eqnarray*}
  \lim_{k\rightarrow\infty}I_{n \times n}+B+B^2+\ldots+B^k = (I_{n \times n}-B)^{-1}.
\end{eqnarray*}

Denoting $\hat{B} \triangleq (I_{n \times n}-B)^{-1}$ we write the distribution of $\mathcal{O}$:
\begin{equation*}
  \mathcal{O} \propto  \hat{B}A \mathcal{I} \sim \mathcal{N}(\hat{B} A \vv, A \hat{B}^T \Sigma_{\vv} \hat{B} A)\,.
\end{equation*}

Recall the definitions of $A$ and $B$ \eqref{AB} and find that:
\begin{equation*}
  \hat{B}A  =(I_{n \times n}-B)^{-1} A = (I_{n \times n}+A W-I_{n \times n})^{-1} A = (A W)^{-1} A = W^{-1}\,.
\end{equation*}
Note that the output sequence can now be characterized more intuitively as:
\begin{equation*}
  \mathcal{O} \propto  W^{-1} \mathcal{I} \sim \mathcal{N}(\uu, W^{-1} \Sigma_{\vv} {W^{-1}}^T)\,.
\end{equation*}

Using known bounds such as multivariate extensions of Chebyshev's bound one can now introduce bounds similar in nature to the box bounds discussed earlier. A better bound can be acquired using the cumulative distribution function for normal distributions; however these can only be estimated numerically. The use of such bounds is:
\begin{equation*}
  p \triangleq \Pr(\|\Ir{k}-\vv\| \leq \delta)\ ,\ \exists \varepsilon \mid \Pr(\|\Or{k}-\uu\| \leq \varepsilon) = p\,,
\end{equation*}
where the connection between $\delta$ and $\varepsilon$ depends on the specific bound used.

So far we have shown that \syncAlg is tolerant towards Gaussian noise, and the output sequence converges to a sequence drawn i.i.d from a normal distribution around the expected solution. 

\subsection{The Asynchronous Setting}

Our proof so far discussed the synchronous case. Next, we extend the analysis of the probabilistic bound for asynchronous case, under mild assumptions. We assume the ``Totally Asynchronous'' model as discussed in \cite{BT97}. We denote $T^i$ as the set of times in which the $i$'th node updates its neighbors with new output values. We further denote $\tau^{i}(t_k)$ as the vector containing the last update times of all received messages at node $i$ at time $t_k$.
\begin{definition}
 {\bf Total Asynchronism}: the sets $T^i$ are infinite and if $\{t_k\}$ is a sequence of elements of $T^i$ that tends to infinity, then $\lim_{k\rightarrow\infty}\tau^{i}_j(t_k) = \infty$ for every $j$.
\end{definition}
This assumption implies that every node is updated infinitely often while allowing maximal flexibility for delay and loss of data.

\cite{BT97} gives a converges proof for the asynchronous case, to an iterative algorithm with constant input values.
The challenge we face is that the input is not constant, but distributed normally. Thus we need to prove that the solution does not infinitely diverge.

For proving that the algorithm output does not diverge, we require some further assumptions. We assume that the update times of each of the nodes are independent of the distribution $\mathcal{I}$. The second assumption we make is that the
input noise distribution is normally distributed, where the noisy reading in different sensors are not correlated.
In other words, the inverse covariance matrix of the input is a diagonal matrix.

To prove that under these assumptions $\mathcal{O}$ is normally distributed we define the set $T$:
\begin{equation*}
  T = \{t \mid \exists i:t \in T^i\}\,.
\end{equation*}
Denote $t_{(i)}=0,\ldots$ as the ordered sequence of times in $T$.
For all $t_{(i)} \in T$ we need to find the distribution of the input at time $t_{(i)}$.
\begin{theorem}\label{theorem:main3}
The infinite input sequence $\Ir{t_{(i)}}$ is distributed normally with $\mathcal{I}$.
\end{theorem}
\begin{proof}
As $\mathcal{I}$ is normally distributed, any marginal distribution computed by a subset of the nodes is normally distributed as well. Now we can induce the requested theorem:

Induction base: By the {\em Total Asynchronism} assumption there exists a time $t_0>0$ at which all nodes have updated their value, and thus the input of each node is now drawn from its marginal distribution. This means that at time $t_0$ the distribution of the entire input vector is (by cartesian product) drawn from $\mathcal{I}$.

Induction hypothesis: At time $t=t_{(i)}$ we assume the input distribution is distributed normally from $\mathcal{I}$.

Induction step: At time $t=t_{(i+1)}$ we separate the nodes into two groups:
  The first group contains nodes that remained constant from the previous round, that is, their a-priori distribution is their marginal distribution drawn from $\mathcal{I}$ by our induction hypothesis.

  The second group contains nodes that are updated with a new input value at this round; this value is distributed again from their marginal distribution.

  Notice that all of these marginal distributions are Gaussian, thus, their cartesian product is Gaussian as well. In addition, any joint covariance between nodes is $0$ due to our assumption. Therefore, we have a cartesian product of all of the nodes distributed as marginal distribution of $\mathcal{I}$ and thus the entire input $\Ir{t_{(i+1)}}$ follows the distribution $\mathcal{I}$.

    Finally we induce that at all $t_{(i)} \in T$: $\Ir{t_{(i)}}$ are drawn from $\mathcal{I}$. Note that unlike the synchronous model - these samples are no longer taken i.i.d from the distribution $\mathcal{I}$.
\qed\end{proof}

Now we can use \theoremref{theorem:main3} and the {\em Total Asynchronism} assumption to draw a similar proof to the synchronous case. Recalling \equationref{eq:recursiveForm}:
\begin{eqnarray*}
  \mathcal{O} \propto \lim_{k\rightarrow\infty}(I_{n \times n}+B+B^2+\ldots+B^{k}) A \mathcal{I}+\lim_{k\rightarrow\infty}B^{k+1} \Or{0}\,.
\end{eqnarray*}
\theoremref{theorem:main3} shows that the input distribution, at all rounds, is drawn from $\mathcal{I}$. This implies that we may apply the Maclaurin series again:
\begin{eqnarray*}
  \lim_{k\rightarrow\infty}I_{n \times n}+B+B^2+\ldots+B^k = (I_{n \times n}-B)^{-1}.
\end{eqnarray*}

We now turn to  \cite[Section 6.3.2]{BT97} where the following claim is proven:
\begin{lemma}\label{claim:contraction}
  For a normalized diagonally dominant matrix $\gamma A$, it holds that
  $x = (I_{n \times n}-\gamma A)x+\gamma b$ asynchronously converges.
\end{lemma}
We now make use of \lemmaref{claim:contraction} under the following notations: $b=0$, which means that our solution converges to $0$. Recalling \equationref{AB}: $B = I_{n \times n}-AW$, where $W$ is diagonally dominant and $A$ is a diagonal matrix, results in $AW$ being diagonally dominant. Combining these results together, we conclude that $B$ is a contraction, thus the linear dynamical system $x=Bx$ converges to zero for any initial $x$:
\begin{eqnarray*}
  \lim_{k\rightarrow\infty}B^{k+1} = 0\,.
\end{eqnarray*}

Regarding convergence speed, the speed at which the initial output factor converges to $0$ and the output distribution converges to a normal distribution, is dependant upon the rate of the least updated node. These settings, though somewhat more restricting than the ones used for synchronous convergence, are all ``minimal'' in the sense that removing any of them creates a system that may diverge. Further discussion about the validity of the
above assumptions is given in \sectionref{sec:discussion}.

\section{Experimental Results}\label{sec:simlation}
\subsection{Box Bound}
\begin{figure}[h!]
\begin{minipage}[b]{0.5\linewidth}
\centering
  \hspace{1.25in} \includegraphics[scale=0.4]{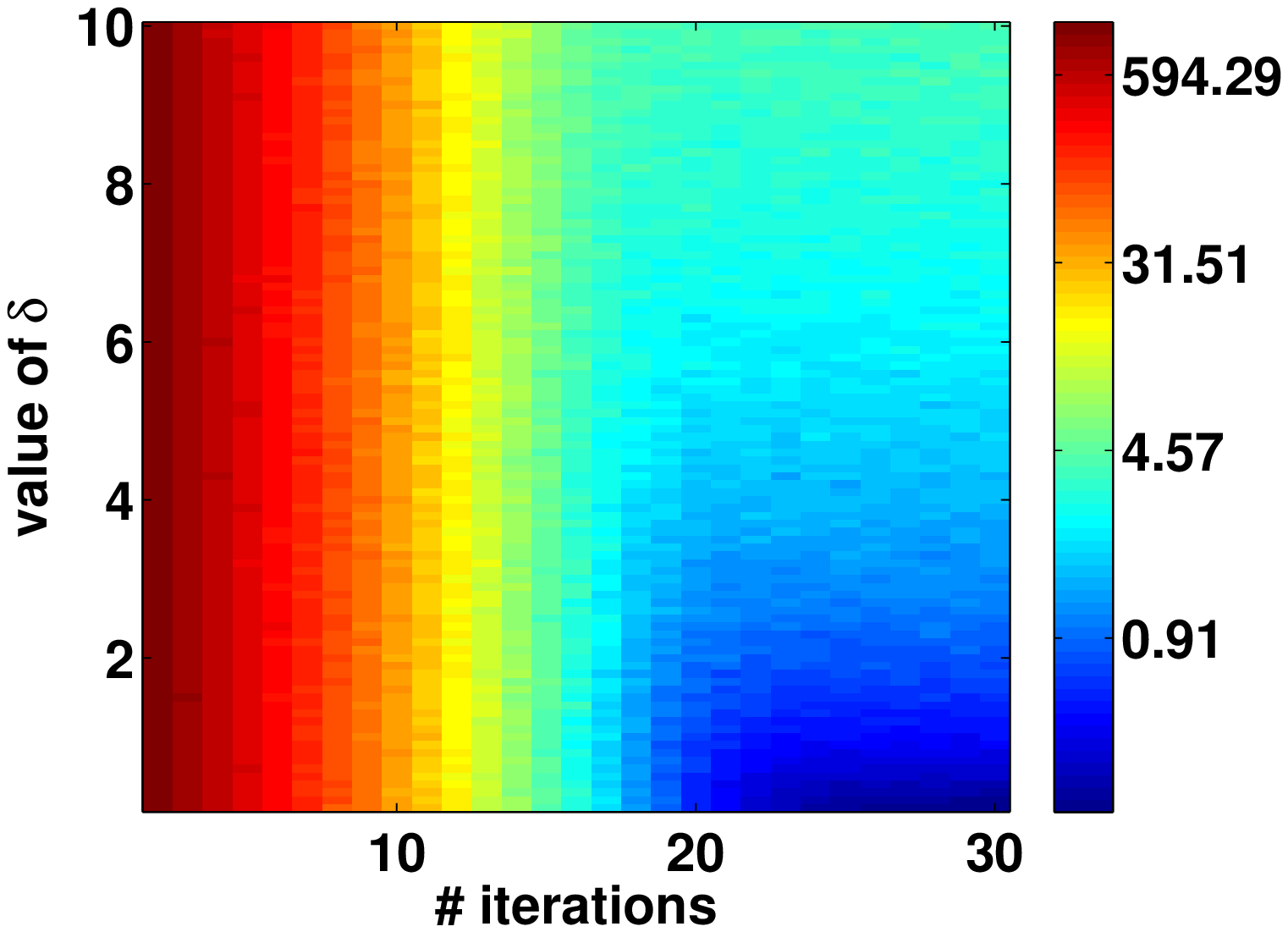}
  \vspace{-0.5cm}
  \caption{Sim. of a Circle graph.}
  \label{figA}
\end{minipage}
\begin{minipage}[b]{0.5\linewidth}
\centering
  \hspace{1.25in} \includegraphics[scale=0.4]{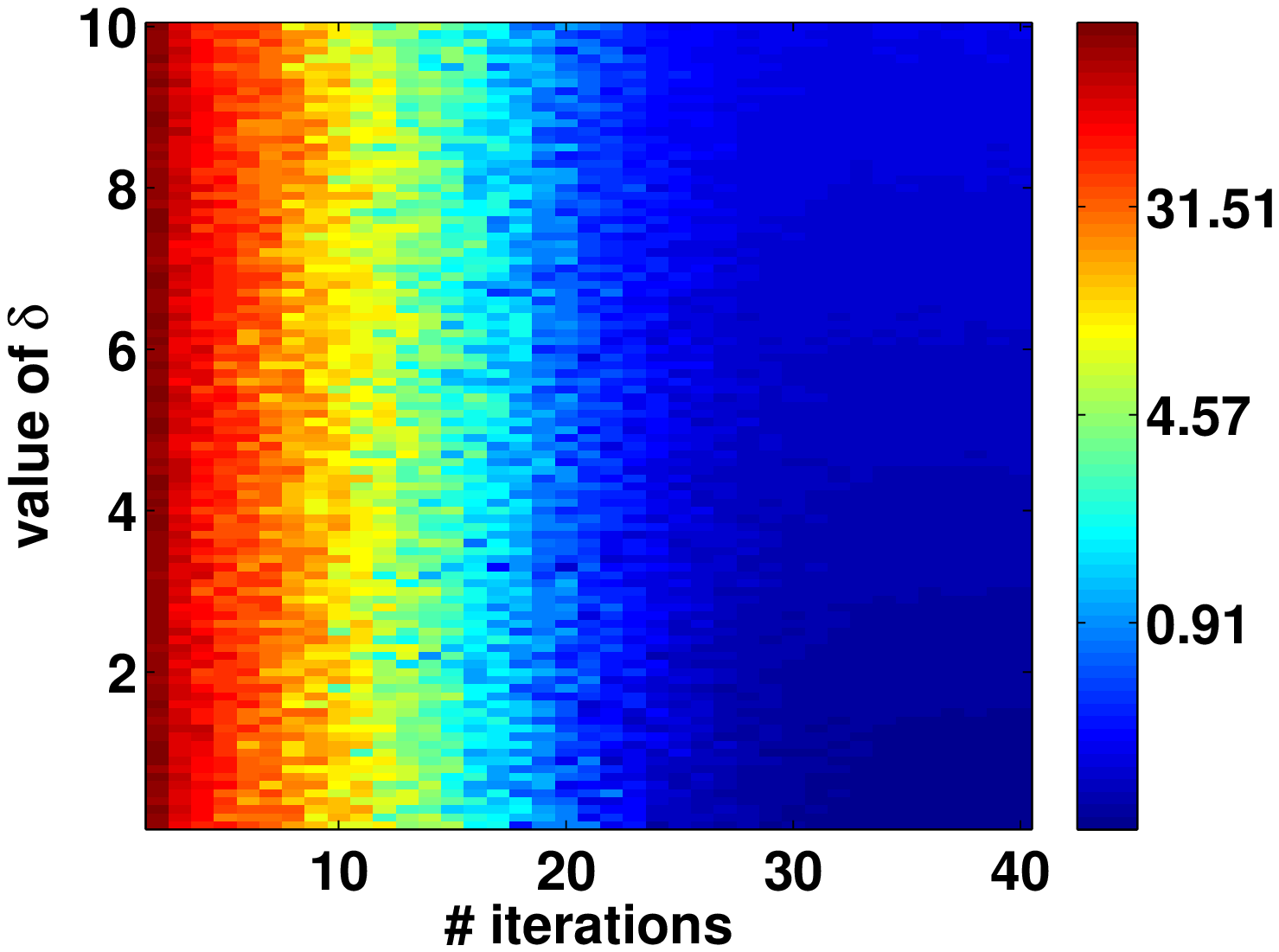}
  \vspace{-0.5cm}
  \caption{Sim. of a Unit-Disc graph.}
  \label{figB}
\end{minipage}
\end{figure}

For illustrating the behavior of our proposed algorithm using the box
bound assumption, we have simulated \syncAlg using two sample topologies of one hundred nodes.
\figureref{figA} depicts a circular topology where each node is
connected to its left and right neighbors. \figureref{figB} shows a
random unit disc graph, where nodes are randomly spread on a
plane, and each node is connected to the nodes that are within a
distance of 1.  The X-axis shows the number of iterations, and the Y-axis shows the value of
$\delta$. Area colors in the heatmap depict the
average of the following procedure: randomly select a vector
$\vv$ and a $\delta$-bounded sequence around $\vv$,
run the simulation for the randomly selected values and return the
$L_\infty$ distance between the last output vector and $\uu$
(calculated as $\uu = W^{-1}\vv$). The heatmap uses a $\log \log$ scale. Both graphs clearly show that as $\delta$ decreases
and the number of iterations increases, the output of \syncAlg
converges to be bounded by a small hypercube around $\uu$.

Note that the unit disc weighted topology matrix is characterized by
$\normi{A}=0.02,\normi{B}=0.97$ while the circle graph
is characterized by $\normi{A}=0.33,\normi{B}=0.66$. As expected,
using unit disc topology requires a larger number of iterations
for convergence (depends on $\normi{B}$). In addition, in the unit
disc topology the value of $\delta$ has a smaller effect on the
convergence, due to the value of $\normi{A}$, which affects the minimal radius around the output.
Since $\normi{A}$ is smaller in the unit disc topology,
increasing $\delta$ does not significantly affect the convergence.

\subsection{Probabilistic Bound}
\begin{figure}[h!]
\begin{minipage}[b]{0.5\linewidth}
\centering
  \hspace{1.25in} \includegraphics[scale=0.3,clip]{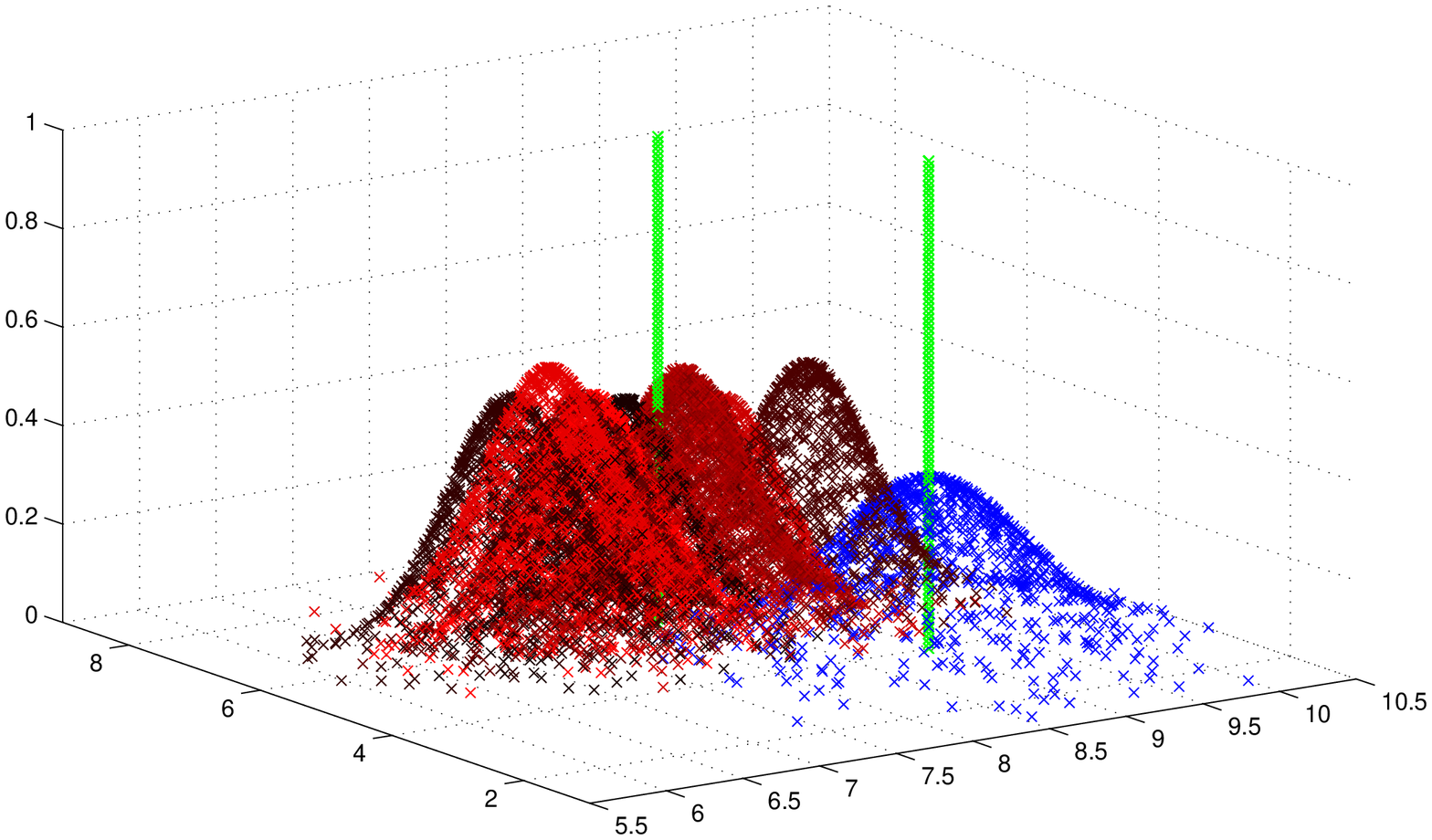}
  \vspace{-0.5cm}
  \caption{Output of several rounds.}
  \label{figC}
\end{minipage}
\begin{minipage}[b]{0.5\linewidth}
\centering
  \hspace{1.25in} \includegraphics[scale=0.3,clip]{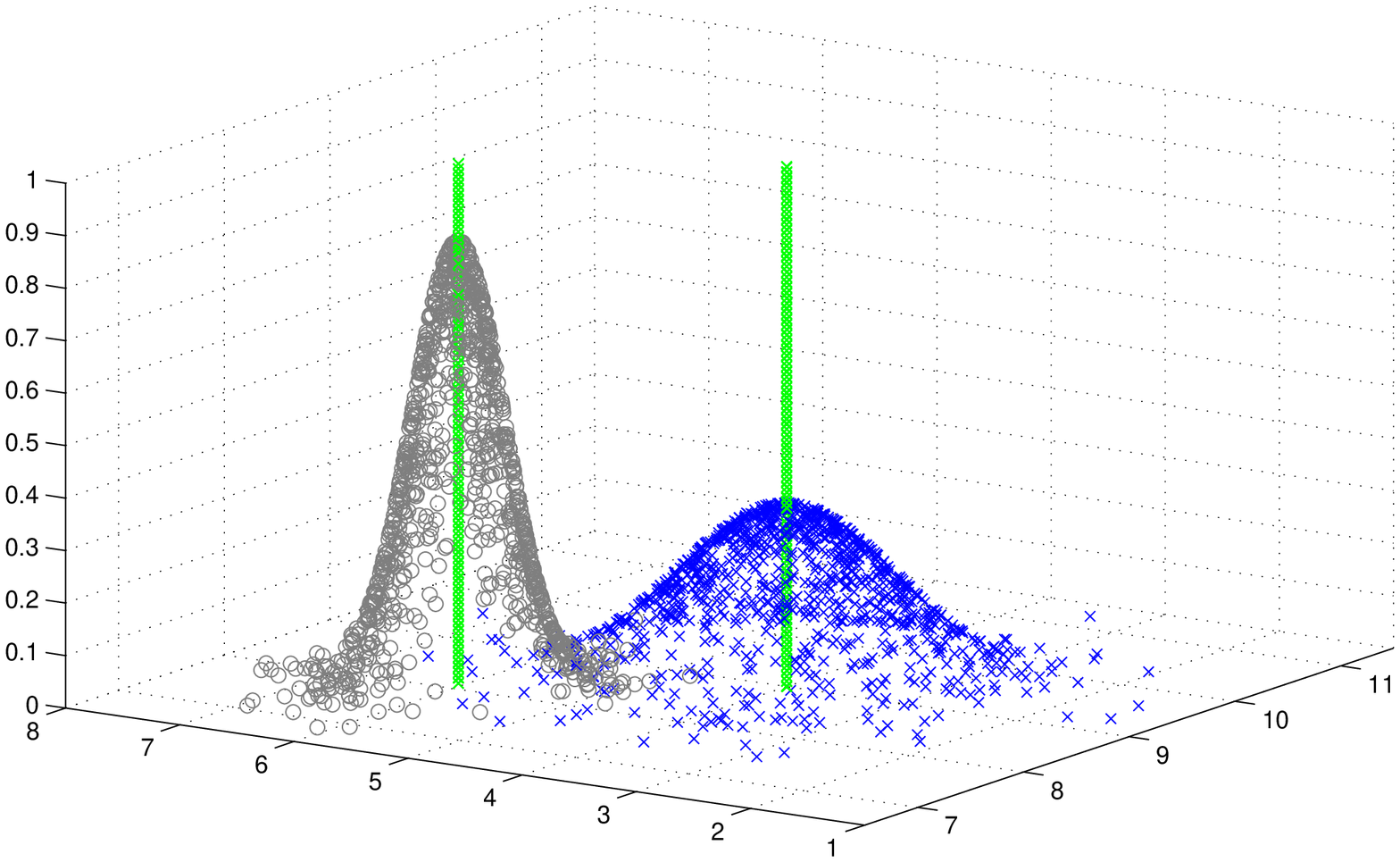}
  \vspace{-0.5cm}
  \caption{Convergence behavior.}
  \label{figD}
\end{minipage}
\end{figure}

We have simulated a run of \syncAlg, when the input noise is i.i.d. normally distributed.
\figureref{figC} plots the behavior of the output distribution at specific rounds. The input distribution of two nodes, printed in blue, is a normal distribution where the two nodes are uncorrelated. The output (in red, growing darker in later rounds) are distributed around some round mean, with the means themselves distributed around $\uu$. It is important to note there is no greedy convergence to $\uu$. This is because each round contains a different noise sample and thus at a specific round our solution can stray away from $\uu$.

\figureref{figD} plots the behavior of the mean value of the input and output to the large system limit. The input
 is to the right (in blue axes) the output is to the left (gray circles). After a large number of rounds, the output is normally distributed around $\uu$, where the plotted $\uu$ (in green) is the result of an EM estimation of the sample mean.

\begin{figure}[h!]
  \centering{
  \includegraphics[scale=0.35,clip]{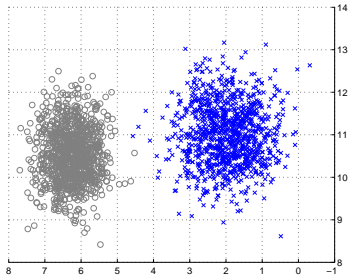}\\
  \caption{Self correlation of input and output.}\label{figE}
  }
\end{figure}

\figureref{figE} plots a similar experiment from a different angle. 
This figure shows the input reading of the nodes (to the right) which is spherically distributed in a circle because the input noise is uncorrelated. The output of the \syncAlg algorithm (to the left) is correlated among nodes, as seen by the elliptic shape of the distribution. The correlation of the output distribution is supported by the theoretical results.

\section{Discussion}\label{sec:discussion}
We have shown that the algorithm \syncAlg is a modification of the
Jacobi iterative method to solve a set of equations $A\xx=\bb$,
where $A$ is given and $\bb$ is dynamically changing but bounded or statistical.
Moreover, \theoremref{theorem:main} is a generalization of
previous analysis of Jacobi's convergence. Our motivation for
\syncAlg originates from the sensor calibration problem, where
sensors need to calibrate their noisy readings. Unlike previous
approaches to this problem, we assume a dynamic system with an
infinite execution of the algorithm. In this setting the readings
of the sensors continuously change. Under the assumption that the
readings' changes are bounded, we have shown that the calibrated
output is bounded as well.

Further application for \syncAlg can be found in any setting where
it is desired to solve $A\xx=\bb$ in a converging and
self-stabilizing manner, while $A$ is given, and $\bb$ may change
slightly from one round to the next. Notice that the analysis
given in
\sectionref{sec:analysis} holds in such a system.

As noted in \remarkref{rem:partofcode} the matrix $A$ is ``part of
the code''. An optional alternative to the current solution is to
compute $A^{-1}$ (the inverted matrix of $A$) beforehand and
include it ``as part of the code''. Thus, each node could locally
solve $\xx = A^{-1}\bb$, and it can be shown that $\xx$ will be
bounded (as long as $\bb$ is bounded). The main problem with such
a solution is the connectivity requirements it incurs. In our
solution, scalar values are sent in the network only between
direct neighbors. The matrix $W$ represents a weighted adjacency
graph. Once inverted, the matrix $W^{-1}$ might not be sparse.
A non-zero entry $w^{-1}_{ij} \in W^{-1}$ means that node $p_i$ needs
to communicate with node $p_j$. This extra communication might cause the algorithm to lose its self-stabilizing properties, as
non-neighboring nodes would require a self-stabilizing overlay
network for their communication.

The assumption of a predefined $A$ is suitable for static networks
in which the communication graph is predetermined. For dynamic
networks, it would be interesting to adjust \syncAlg to discover
the connectivity of the network, inferring the optimal weights dynamically.
We assume that after the weights are calculated, the topology of the
sensor network remains stable, thus the convergence analysis of \sectionref{sec:analysis} should hold.

\subsection{Necessity of Convergence Conditions for the Probabilistic Bound}
Assuming that the sensor input readings are normally distributed due to measurement noise, the resulting output is also normally distributed. As discussed in \sectionref{sec:Probabilistic}, this is not a sufficient condition for the convergence of \syncAlg in a finite number of rounds.

We give the following example to demonstrate the claim above. Assume the input reading at round $1$ is larger than $\vv$, and at each of the following rounds the input reading of each node is larger than at the previous round. As the input is normally distributed, for every $n>0$ there is a small but positive chance for such a sequence of length $n$. Input sequence behaving as described would lead \syncAlg to a growing output sequence for $n$ rounds, meaning for the duration of these rounds \syncAlg diverges. Thus it is obvious the mean can be infinitely far from $\uu$ as long as $n$ is finite. However, when $n$ is large, the probably of such an anomaly is diminishing to zero.

The convergence results for the asynchronous case had a minimal set of assumptions. One of them is that the input noise distribution must be uncorrelated among nodes. That is because if the input noise distribution is correlated, all nodes
must sample the input concurrently, an assumption which is not valid under the asynchronous model. 

\subsection{Relation to Perturbation Theory}
A large amount of research focused on the problem of solving $A\xx=\bb$ when $A$ and
$\bb$ are not exactly known. That is, let $\hat{A}=A+\delta A$ and
$\hat{\bb} = \bb+\delta \bb$, and consider the equation
$\hat{A}\hat{\xx}=\hat{\bb}$; what can be said about $\xx$ in
relation to $\hat{\xx}$?

Our setting is ``easier'' in one sense, and ``harder'' in a
different sense. In our setting $A$ is known, \ie $\delta A
= 0$. However, $\hat{\bb}$ is not well defined. That is, the input
vector - which is described by $\hat{\bb}$ - changes over time.
When solving $\hat{A}\hat{\xx}=\hat{\bb}$ it is assumed that there
is some $\bb$ that is \emph{constant} but it was measured with an
error. In our case, $\bb$ is not constant as it changes over time,
while it represents the measurements correctly.

As a future research, it would be interesting to consider the
implications of adding inaccuracy to the measurements. The vast
body of knowledge regarding perturbation theory would definitely
aid in this extension to our model.

\subsection{Relation to Convex Optimization}
Many practical optimization problems are given in the quadratic
form $f(x) = 1/2\xx A\xx - \bb^T\xx$, where the task is to compute
$\min_{\xx} f(\xx)$ distributively over a communication network. A
survey showing several applications can be found in \cite{PPNA08}.
Example applications are monitoring, distributed computation of
trust and ranking of nodes and data items.

A standard way for solving $\min_{\xx} f(\xx)$ is by computing the
derivative and comparing it to zero to get the global optimum.
When the matrix $A$ is symmetric, $f'(\xx) = A\xx -\bb = 0$, and we get
a linear system of equations $A\xx = \bb$. In other words, the convex
optimization problem is reduced into a solution of a linear system
of equations.

Interior point methods \cite[Ch. 11]{BV04} solve linear programming
problems by applying Newton method iteratively. Each computation
of the Newton step involves a solution of a linear systems of
equations. An area of future work is to examine the applications
of our self-stabilizing algorithm to these methods. The
difficulties arise from the fact that the matrix A needs to be
recomputed between iterations, so nodes need to be synchronized
and aware of the current iteration taking place.

\section*{Acknowledgements}
The authors would like to thank Golan Pundak for
assisting with the simulations, and the anonymous reviewers of SSS-08 conference
for their helpful comments.

\bibliographystyle{plain}
\bibliography{TCS09}




\end{document}